\documentclass{article}

\usepackage[utf8]{inputenc}
\usepackage[english]{babel}

\usepackage{graphicx}
\usepackage{amssymb}
\usepackage{amsmath}
\usepackage{tikz-cd}
\usepackage{enumerate}
\usepackage{verbatim}
\usepackage{tikz}
\usepackage{ytableau}
\usepackage{hyperref}

\newtheorem{theorem}{Theorem}[section]
\newtheorem{corollary}[theorem]{Corollary}
\newtheorem{lemma}[theorem]{Lemma}
\newtheorem{proposition}[theorem]{Proposition}
\newtheorem{conjecture}[theorem]{Conjecture}

\newcommand{\End}{\text{End}}

\newcommand{\Ind}{\text{Ind}}
\newcommand{\Res}{\text{Res}}
\newcommand{\Sym}{\text{Sym}}

\let\emptyset\varnothing

\begin{document}

\title{Representations of monomial matrices and restriction from $GL_n$ to $S_n$}
\author{{Nate Harman}}


\maketitle

\vspace{-.9cm}
\begin{abstract}

We outline an approach to understanding restrictions of polynomial representations of $GL_n(\mathbb{C})$ to $S_n$ by first restricting to $T \rtimes S_n$, the subgroup of $n \times n$ monomial matrices. Using this approach we give a combinatorial interpretation for the decomposition of a tensor product of symmetric powers of the defining representation.

\end{abstract}

\begin{section}{Introduction and notation}

The symmetric group $S_n$ sits inside the general linear group $GL_n = GL_n(\mathbb{C})$ as the set of permutation matrices.  Suppose $\lambda$ is a partition of size at most $n$, we can view it as a dominant polynomial weight for the diagonal torus $T \subset GL_n$. Let $W(\lambda)$ denote the corresponding irreducible polynomial $GL_n$ representation of highest weight $\lambda$.  If $\mu$ is a partition of $n$ one may consider the following question:

\medskip

\noindent \textbf{Problem:} Give a positive combinatorial interpretation for the multiplicity of the irreducible Specht module $S^\mu$ inside the restriction of $W(\lambda)$ to $S_n$.

\medskip

In terms of symmetric functions this problem is closely related to the problem of understanding the decompositions of certain plethysms and inner plethysms  into Schur functions, a notoriously mysterious problem in combinatorial representation theory.

\medskip

Recall that $W(\lambda)$ decomposes as a direct sum of $T$-weight spaces $W(\lambda)_\nu$.  In general, the action of $S_n$ does not preserve these weight spaces, however it does permute them according to the natural action of $S_n$ on the weights of $T$. 

In particular if we let $\overline{W(\lambda)_\nu}$ denote the direct sum of weight spaces for weights that are $S_n$-conjugate to $\nu$ then this carries a natural action of $S_n$. Hence we can further refine the above problem to the following problem:

\medskip

\noindent \textbf{Problem:} Give a positive combinatorial interpretation for the multiplicity of the irreducible Specht module $S^\mu$ inside $\overline{W(\lambda)_\nu}$.

\medskip

If we let $T \rtimes S_n = N(T) = \mathbb{C}^\times \wr S_n$ denote the subgroup of monomial matrices (matrices with a single non-zero entry in each row and column) this problem is essentially equivalent to the following two steps:

\begin{enumerate}

\item Understand how $W(\lambda)$ decomposes as a $T \rtimes S_n$ representation.

\item Understand how irreducible $T \rtimes S_n$ representations decompose when restricted to $S_n$.

\end{enumerate}

In this paper we give some partial progress in carrying out this program. The rest of the paper is organized as follows:

\medskip

 In Section \ref{primer} we go over some generalities about algebraic representations of $T \rtimes S_n$, define the class of polynomial representations, characterize the irreducible algebraic and polynomial representations of $T \rtimes S_n$, and explain completely how to do step 2 above. This part will be mostly expository and will largely follow from standard facts about representations of wreath products.
 
 \medskip
 
 In Section \ref{combo} we define the class of weighted permutation modules of $T \rtimes S_n$ and describe some combinatiorics involving them. We then use them to give the main result of this paper, a combinatorial interpretation for the decomposition of the restriction to $T \rtimes S_n$ of $$Sym^{k_1}(V) \otimes Sym^{k_2}(V) \otimes \dots \otimes Sym^{k_m}(V)$$ where $V$ denotes the defining $n$-dimensional representation of $GL_n$. As a corollary this gives a combinatorial description of the restriction to $S_n$ as well.
 
 \medskip
 
In Section \ref{future} we suggest some future directions of research by outlining a version of Schur-Weyl duality for $T \rtimes S_n$, explain a connection of this work to Foulkes conjecture about plythysm, and give some explicit calculations of restrictions of low degree $GL_n$ representations to $T \rtimes S_n$.

\end{section}

\begin{section}*{Acknowledgements}
Thanks to Benson Farb and the UChicago geometry, topology, and representation theory working group. This work was partially supported by NSF postdoctoral fellowship award 1703942.

\end{section}

\begin{section}{Representation theory of $T \rtimes S_n$ - a primer}\label{primer}

In this section we will give an overview of the algebraic and polynomial representation theory of $T \rtimes S_n$ and the related combinatorics for working with these representations.

\begin{subsection}{Algebraic representations of $T \rtimes S_n$}

The group $T \rtimes S_n$ is a wreath product of $S_n$ with $\mathbb{C}^\times$, and as such its (algebraic) representations are easy to describe. If $G$ is a finite group the representation theory of $G \wr S_n$ is standard material and can be found in a number of basic representation theory texts (in particular it is explained in detail in \cite{Zelevinsky}) and the theory goes through mutatis mutandis if instead $G$ is a reductive algebraic group over $\mathbb{C}$ and we are considering algebraic representations of $G \wr S_n$.

\subsubsection{Irreducible algebraic representations}

The simplest irreducible algebraic representations are the weighted Specht modules $S^{\lambda,k}$ where as an $S_n$ representation this is just the usual Specht module $S^\lambda$, and each copy of $\mathbb{C}^\times$ just acts by scalars via the character $z \to z^k$.

If $\lambda^1, \lambda^2, \dots, \lambda^m$ are partitions of total size $n$, and $k_1 < k_2 < \dots < k_m$ are distinct integers then the induced representation 
\begin{equation}\label{algirr}
\text{Ind}_{(\mathbb{C}^\times \wr S_{|\lambda^1|}) \times (\mathbb{C}^\times \wr S_{|\lambda^2|}) \times \dots (\mathbb{C}^\times \wr S_{|\lambda^m|})}^{ \mathbb{C}^\times \wr S_{n} }(S^{\lambda^1,k_1} \otimes S^{\lambda^2,k_2} \otimes \dots \otimes S^{\lambda^m,k_m})
\end{equation}
is irreducible and moreover every irreducible algebraic representation of $T \rtimes S_n$ is obtained this way. Hence we see the irreducible algebraic representations of $T \rtimes S_n$ can be naturally labeled by collections of partitions of total size $n$ indexed by the integers .

\subsubsection{Induction of representations}

Moreover, it is easy to describe the induction of two weighted Specht modules of the same $\mathbb{C}^\times$ weight. We just get a direct sum of weighted Specht modules of same $\mathbb{C}^\times$  weight with multiplicities coming from the usual Littlewood-Richardson rule. More precisely: 
$$\text{Ind}_{(\mathbb{C}^\times \wr S_{|\lambda|}) \times (\mathbb{C}^\times \wr S_{|\mu|})} ^ {\mathbb{C}^\times \wr S_{|\lambda| + |\mu|} }(S^{\lambda,k} \otimes S^{\mu,k}) = \bigoplus_\nu  c_{\lambda, \mu}^\nu S^{\nu, k}$$
It's then clear how to decompose the induction of two arbitrary irreducibles, you just run the Littlewood-Richardson rule separately on each $\mathbb{C}^\times$-weight. In other words the Grothendieck ring of $\bigoplus_n \text{Rep}(\mathbb{C}^\times \wr S_{n})$ with monoidal structure coming from induction is just a direct sum indexed by $\mathbb{Z}$ of copies of the ring of symmetric functions $\Lambda$.

\subsubsection{Restriction to $S_n$}

Finally, we will close out this section by noting that because irreducible representations of $T \rtimes S_n$ are just induced up from products of weighted Specht modules, we can easily describe the restriction of arbitrary irreducible $T \rtimes S_n$ representations to $S_n$. The irreducible representation 

$$\text{Ind}_{(\mathbb{C}^\times \wr S_{|\lambda^1|}) \times (\mathbb{C}^\times \wr S_{|\lambda^2|}) \times \dots (\mathbb{C}^\times \wr S_{|\lambda^m|})}^{ \mathbb{C}^\times \wr S_{n} }(S^{\lambda^1,k_1} \otimes S^{\lambda^2,k_2} \otimes \dots \otimes S^{\lambda^m,k_m})$$
just becomes
$$\text{Ind}_{S_{|\lambda^1|} \times S_{|\lambda^2|} \times \dots S_{|\lambda^m|}}^{S_{n} }(S^{\lambda^1} \otimes S^{\lambda^2} \otimes \dots \otimes S^{\lambda^m})$$
which of course can be decomposed using the Littlewood-Richardson rule.  Explicitly in terms of symmetric functions we have the following corollary.

\begin{corollary}\label{restmult}
The multiplicity of the Specht module $S^\mu$ inside the restriction to $S_n$ of $$\text{Ind}_{(\mathbb{C}^\times \wr S_{|\lambda^1|}) \times (\mathbb{C}^\times \wr S_{|\lambda^2|}) \times \dots (\mathbb{C}^\times \wr S_{|\lambda^m|})}^{ \mathbb{C}^\times \wr S_{n} }(S^{\lambda^1,k_1} \otimes S^{\lambda^2,k_2} \otimes \dots \otimes S^{\lambda^m,k_m})$$ is equal to the multiplicity of the Schur function $s_\mu$ in the product $s_{\lambda^1}s_{\lambda^2}\dots s_{\lambda^m}$

\end{corollary}

As suggested in the introduction we see that the difficulty in understanding the $S_n$ action on $\overline{W(\lambda)_\nu}$ lies in understanding it as a $T \rtimes S_n$ representation, and that then restricting to $S_n$ is straightforward.

\end{subsection}

\begin{subsection}{Polynomial representations of $T \rtimes S_n$}

We say that an algebraic representation of $T \rtimes S_n$ is \emph{weakly polynomial} if the $T$-weights that occur are all polynomial weights. The weighted Specht module $S^{\lambda, k}$ is weakly polynomial if and only if $k \ge 0$, and more generally an irreducible representation of the form (\ref{algirr}) is weakly polynomial if and only if $k_i \ge 0$ for all $i$.  The class of weakly polynomial representations is preserved by induction and the corresponding Grothendieck ring of $\bigoplus_n \text{Rep}^{wp}(\mathbb{C}^\times \wr S_{n})$ is again a direct sum of copies of $\Lambda$, this time indexed by $\mathbb{Z}_{\ge 0}$.

Let $V$ denote the standard $n$-dimensional representation of $T \rtimes S_n$ obtained by restricting the defining representation of $GL_n$ to $T \rtimes S^n$. As a representation of $S_n$ this is just the standard $n$-dimensional permutation representation, and $T$ acts on the permutation basis with weights $(1,0,0,\dots, 0)$, $(0,1,0,0,\dots 0)$, $\dots$, and $(0,0,\dots, 0, 1)$. In the notation of (\ref{algirr}) we have that $$V \cong \text{Ind}_{(\mathbb{C}^\times \wr S_{n-1}) \times (\mathbb{C}^\times \wr S_{1})} ^ {\mathbb{C}^\times \wr S_{n} }(S^{(n-1),0} \otimes S^{(1),1}).$$

We say that a representation of $T \rtimes S_n$ is \emph{polynomial of degree $d$} if it is a direct summand of copies of $V^{\otimes d}$, and say that a representation is \emph{polynomial} if it is a direct sum of polynomial representations of various degrees.  In particular, a representation of $T \rtimes S_n$ is polynomial if and only if it is a direct summand of the restriction of a polynomial representation of $GL_n$. 

\medskip

\noindent \textbf{Warning:} We'll note that unlike the classes of algebraic and weakly polynomial representations, induction from $(\mathbb{C}^\times \wr S_{n}) \times (\mathbb{C}^\times \wr S_{m})$ to $\mathbb{C}^\times \wr S_{n+m}$ does not preserve the class of polynomial representations.  However it is obvious by the definition that polynomial representations are preserved by the internal tensor product of $T \rtimes S_n$ representations.

\subsubsection{Irreducible polynomial representations}

The following proposition characterizes which irreducible weakly polynomial representations are polynomial.

\begin{proposition}\label{polirr}
Given partitions $\lambda^1, \lambda^2, \dots, \lambda^m$ of total size $n$, and distinct nonnegative integers $0 \le k_1 < k_2 < \dots < k_m$ then the irreducible weakly polynomial representation

$$\text{Ind}_{(\mathbb{C}^\times \wr S_{|\lambda^1|}) \times (\mathbb{C}^\times \wr S_{|\lambda^2|}) \times \dots (\mathbb{C}^\times \wr S_{|\lambda^m|})}^{ \mathbb{C}^\times \wr S_{n} }(S^{\lambda^1,k_1} \otimes S^{\lambda^2,k_2} \otimes \dots \otimes S^{\lambda^m,k_m})$$
is polynomial (of degree $k_1+k_2+\dots k_m$) if and only if either $k_1 \ne 0$ or $k_1 = 0$ and $\lambda^1$ is a one-row partition.

\end{proposition}

Before proving Proposition \ref{polirr} it will be convenient to introduce a new labeling for these representations.  If $\lambda^1, \lambda^2, \dots, \lambda^\ell$ is a sequence of (possibly empty) partitions of total size $m \le n$ then define $$V^{\lambda^1, \lambda^2, \dots, \lambda^\ell} := \text{Ind}_{(\mathbb{C}^\times \wr S_{n-m}) \times (\mathbb{C}^\times \wr S_{|\lambda^1|}) \times \dots (\mathbb{C}^\times \wr S_{|\lambda^\ell|})}^{ \mathbb{C}^\times \wr S_{n} }(S^{(n-m),0} \otimes S^{\lambda^1,1} \otimes \dots \otimes S^{\lambda^\ell,\ell})$$

For example in this notation the standard $n$-dimensional representation $V$ from the previous section is now denoted by $V^{(1)}$, and if instead we took the $n$-dimensional permutation representation for $S_n$ and had $T$ act on the permutation basis with weights $(2,0,0,\dots, 0)$, $(0,2,0,0,\dots 0)$, $\dots$, and $(0,0,\dots, 0, 2)$ this representation would be labeled $V^{\emptyset, (1)}$.  Under this labeling the degree of $V^{\lambda^1, \lambda^2, \dots, \lambda^\ell}$ is given by $|\lambda^1| + 2|\lambda^2| + 3|\lambda^3| + \dots + \ell | \lambda^\ell|.$

We'll note that adding empty partitions to the end of the sequence $\lambda^1, \lambda^2, \dots, \lambda^\ell$ doesn't change the corresponding representation.  It will occasionally be useful to think of this as being a labeling by infinite sequences of partitions where all but finitely many are the empty partition.  In particular, the trivial representation just corresponds to the sequence where all the partitions are empty.

At times we will want to use this notation to describe representations of $\mathbb{C}^\times \wr S_{n}$ for different values of $n$, in which case we will add subscripts $V_{n}^{\lambda^1, \lambda^2, \dots, \lambda^\ell}$ to avoid ambiguity.

\medskip

Since by definition all irreducible polynomial representations are direct summands of tensor powers of $V^{(1)}$ it will be useful to see how tensoring with $V^{(1)}$ looks under this labeling.  The following lemma gives such a description:

\begin{lemma}\label{tensorV} \hspace{1cm} \\

 \noindent If $|\lambda^1| + |\lambda^2| + \dots + |\lambda^\ell| = n$ then  $$V^{(1)} \otimes V^{\lambda^1, \lambda^2, \dots, \lambda^\ell} = \bigoplus_{\substack{{\lambda^i}' = \lambda^i - \square \\ {\lambda^{i+1}}' = \lambda^{i+1} + \square}} V^{\lambda^1, \dots ,{\lambda^i}' , {\lambda^{i+1}}',  \dots, \lambda^\ell}$$
and if $|\lambda^1| + |\lambda^2| + \dots + |\lambda^\ell| < n$

$$V^{(1)} \otimes V^{\lambda^1, \lambda^2, \dots, \lambda^\ell} = \bigoplus_{\substack{{\lambda^i}' = \lambda^i - \square \\ {\lambda^{i+1}}' = \lambda^{i+1} + \square}} V^{\lambda^1, \dots ,{\lambda^i}' , {\lambda^{i+1}}',  \dots, \lambda^\ell} \oplus \bigoplus_{{\lambda^1}' = \lambda^1 + \square} V^{{\lambda^1}', \lambda^2,\dots \lambda^\ell}$$
where the sums are over all ways of removing/adding a box to the corresponding Young diagrams.
\end{lemma}

In terms of our labeling by sequences of partitions this just says that tensoring with $V^{(1)}$ corresponds to the following process:  If we start with one sequence of partitions we create new ones in all possible ways by removing one box from one partition and adding it to the next partition in the sequence, or if there aren't already $n$ total boxes we can add one box to the first partition.

We'll mostly be interested in the case where $n$ is strictly larger than the degree so we are always allowed to add a box to the first partition, but we've included the general case here for completion. In terms of the Young diagrams involved here is a picture of the first few rows of the Bratteli diagram for tensoring with $V^{(1)}$ when $n$ is at least $3$ (when $n$ is less than $3$ just delete all sequences of Young diagrams with more than $n$ boxes):

\medskip

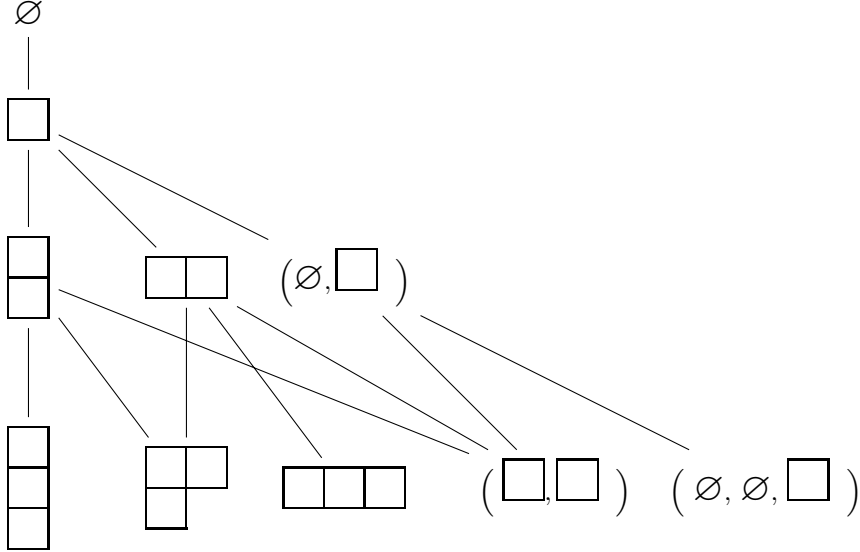
\begin{tikzpicture}[scale=.7]
  \node (empty) at (0,0) {\scalebox{1.5}{$\emptyset$}};
  \node (1) at (0,-2) {$\ydiagram{1}$};
  \node (a) at (0,-5) {$\ydiagram{1,1}$};
  \node (b) at (3,-5) {$\ydiagram{2}$};
  \node (c) at (6,-5) {$\Big{(} \scalebox{1.5}{$\emptyset$},  \ydiagram{1}$ \ \Big{  )}};
  \node (d) at (0,-9) {$\ydiagram{1,1,1}$};
   \node (e) at (3,-9) {$\ydiagram{2,1}$};
 \node (f) at (6,-9) {$\ydiagram{3}$};
 \node (g) at (10,-9) {$\Big{(} \ \ydiagram{1},  \ydiagram{1}$ \ \Big{  )}};
  \node (h) at (14,-9) {$\Big{(}$ \scalebox{1.5}{$\emptyset$}, \scalebox{1.5}{$\emptyset$},  $\ydiagram{1}$ \ \Big{  )}};
  
  \draw (empty)--(1)--(a)--(d)
  (1)--(b)--(e)
  (1) -- (c)--(g)
  (a)--(e)
  (b)--(f)
  (a)--(g)--(b) 
 (c)--(h);
\end{tikzpicture}

\medskip

\noindent \textbf{Proof of Lemma \ref{tensorV}:} $V^{(1)}$ is induced from $(\mathbb{C}^\times \wr S_{n-1}) \times (\mathbb{C}^\times \wr S_{1})$ so in order to tensor with it we may use the push-pull formula:

$$\Ind(U) \otimes W \cong \Ind(U \otimes \Res(W))$$
 
So to calulate this first we need to restrict $V^{\lambda^1, \lambda^2, \dots, \lambda^\ell}$ to $(\mathbb{C}^\times \wr S_{n-1}) \times (\mathbb{C}^\times \wr S_{1})$. We already saw that induction between algebraic representations was controlled by the Littlewood-Richardson rule, hence by Frobenius reciprocity it follows that if $|\lambda^1| + |\lambda^2| + \dots + |\lambda^\ell| = n$ then

$$\Res (V_n^{\lambda^1, \lambda^2, \dots, \lambda^\ell}) = \bigoplus_{{\lambda^i}' = \lambda^i - \square} V_{n-1}^{\lambda^1, \dots {\lambda^i}' \dots, \lambda^\ell} \otimes S^{(1),i}$$

and if $|\lambda^1| + |\lambda^2| + \dots + |\lambda^\ell| < n$ then 

$$\Res  (V_n^{\lambda^1, \lambda^2, \dots, \lambda^\ell}) = \bigoplus_{{\lambda^i}' = \lambda^i - \square} V_{n-1}^{\lambda^1, \dots {\lambda^i}' \dots, \lambda^\ell} \otimes S^{(1),i} \hspace{.4cm} \oplus \hspace{.4cm}  V_{n-1}^{\lambda^1, \lambda^2, \dots, \lambda^\ell} \otimes S^{(1),0}$$

Next we need to tensor these factors with the $(\mathbb{C}^\times \wr S_{n-1}) \times (\mathbb{C}^\times \wr S_{1})$ representation that is the trivial representation of $(\mathbb{C}^\times \wr S_{n-1})$ tensored with $S^{(1),1}$. Of course tensoring with the trivial representation does nothing, and the second factor is just a character of $\mathbb{C}^\times$ so we see this just sends $V_{n-1}^{\mu^1, \dots, \mu^\ell} \otimes S^{(1),i}$ to $V_{n-1}^{\mu^1, \dots, \mu^\ell} \otimes S^{(1),i+1}$.

Finally, we need to induce back up to $(\mathbb{C}^\times \wr S_{n-1})$.  But we know by the Littlewood-Richardson rule that $$\Ind V_{n-1}^{\mu^1, \dots, \mu^\ell} \otimes S^{(1),i+1} = \bigoplus_{{\mu^{i+1}}' = \mu^{i+1} + \square} V_n^{\mu^1, \dots {\mu^{i+1}}', \dots, \mu^\ell}$$
Putting this all together gives the desired formula. \hfill $\square$

\medskip

From here Proposition \ref{polirr} follows immediately: By induction every irreducible subrepresentation of $(V^{(1)})^{\otimes d}$ is of the desired form, and conversely it's clear that any sequence of partitions can be obtained from the empty sequence by the process described in Lemma \ref{tensorV} of adding boxes to the first partition and moving them down the sequence arbitrarily (moreover the multiplicity is given by the number of downward walks on the Bratteli diagram from the empty sequence to the desired sequence of partitions).

\medskip

\noindent \textbf{Example:} If $n$ is at least $2$ Lemma \ref{tensorV} gives us the decomposition of the tensor square of the defining representation:

$$V^{(1)} \otimes V^{(1)} \cong V^{(1,1)} \oplus V^{(2)} \oplus V^{\emptyset, (1)}$$
In terms of $GL_n$ representations this is familiar, if we take the usual decomposition of $GL_n$ representations
$$V \otimes V = \Lambda^2(V) \oplus \Sym^2(V)$$
then $\Lambda^2(V)$ restricts to $V^{(1,1)}$, and if we think of $\Sym^2(V)$ as degree $2$ polynomials of standard basis vectors $x_1, x_2, \dots x_n$ then $V^{\emptyset, (1)}$ corresponds to the span of the polynomials $x_i^2$ and $V^{(2)}$ corresponds to the span of the polynomials $x_ix_j$ with $i \ne j$.

\end{subsection}

\end{section}

\begin{section}{Permutation modules, symmetric powers, and multiset combinatorics}\label{combo}

\begin{subsection}{Weighted permutation modules}

Recall for the symmetric group $S_n$ and a partition (or composition) $\lambda =(\lambda_1,\lambda_2, \dots, \lambda_\ell)$, the permutation module $M(\lambda)$ is defined as $$M(\lambda) := \text{Ind}_{S_{\lambda_1} \times  S_{\lambda_2} \times \dots  S_{\lambda_\ell}}^{S_n}(\mathbf{1}) = \text{Ind}_{S_{\lambda_1} \times  S_{\lambda_2} \times \dots  S_{\lambda_\ell}}^{S_n}(S^{(\lambda_1)} \otimes S^{(\lambda_2)}\otimes \dots \otimes S^{(\lambda_\ell)}  )$$
which combinatorially just means $M(\lambda)$ is the linearized representation of the action of $S_n$ on ordered set partitions of $\{1,2,\dots,n\}$ into sets $A_1, A_2, \dots A_\ell$ with $|A_i| = \lambda_i$ for each $i$.

For algebraic representations of $T \rtimes S_n$ there are natural analogs of these modules where we allow the copies of $\mathbb{C}^\times$ to act by characters and then induce up. Explicitly we have the \emph{weighted permutation modules} $M(\lambda, \mathbf{k})$ where $\lambda = (\lambda_1,\lambda_2, \dots, \lambda_\ell)$ is a composition of $n$ and $\mathbf{k} = (k_1, k_2, \dots, k_\ell)$ is a list of $\mathbb{C}^\times$ weights of the same length
$$M(\lambda, \mathbf{k}) := \text{Ind}_{(\mathbb{C}^\times \wr S_{|\lambda_1|}) \times (\mathbb{C}^\times \wr S_{|\lambda_2|}) \times \dots (\mathbb{C}^\times \wr S_{|\lambda_\ell|})}^{ \mathbb{C}^\times \wr S_{n} }(S^{(\lambda_1),k_1} \otimes S^{(\lambda_2),k_2} \otimes \dots \otimes S^{(\lambda_\ell),k_\ell})$$
Note that we can rearrange the parts of the composition or the order of the $k_i$'s, so long as we permute the other accordingly. In particular, at times it will be convenient to reorder things so that $\lambda$ is a partition and other times it will be convenient to order the $k_i$'s in increasing order. Just note you can't in general do both simultaneously.

It's then easy to see that a weighted permutation module is polynomial if and only if each $k_i$ is non-negative and at most one of them is equal to zero.  In this case it will be convenient to reindex as before by lists of partitions $\lambda^1, \lambda^2, \dots \lambda^j$ of total size at most $n$. If we denote $\lambda^i = (\lambda^i_1, \lambda^i_2,\dots, \lambda^i_{\ell_i})$ and $|\lambda^1|+ \dots + |\lambda^j| = m$ then define
$$\tilde{M}(\lambda^1, \lambda^2,\dots, \lambda^k) := M((n-m, \lambda^1_1, \lambda^1_2, \dots, \lambda^j_{\ell_j}), (0, 1,1,\dots, 2,2,\dots, j,j))$$
where there are $\ell_1$ $1$'s, $\ell_2$ $2$'s, and so on. 

 In other words, $\tilde{M}(\lambda^1, \lambda^2,\dots, \lambda^k)$ is just the representation induced from the representation of $(\mathbb{C}^\times \wr S_{n-m})\times (\mathbb{C}^\times \wr S_{|\lambda^1|}) \times \dots (\mathbb{C}^\times \wr S_{|\lambda^j|})$ where the first factor acts trivially and in the $\mathbb{C}^\times \wr S_{|\lambda^i|}$ factor we take a copy of the $S_n$ permutation module $M(\lambda^i)$ with the copies of $\mathbb{C}^\times$ scaling by the character $z \to z^i$.
 
\medskip

\noindent \textbf{Examples:} The defining representation $V^{(1)}$ is isomorphic to the permutation module $\tilde{M}((1))$.  If $n \ge 2$ the tensor square $V^{(1)} \otimes V^{(1)}$ decomposes as $\tilde{M}((1,1)) \oplus \tilde{M}(\emptyset, (1))$. In terms of the standard basis vectors $x_1, x_2, \dots x_n$ of $V^{(1)}$, $\tilde{M}((1,1))$ is the span of all vectors of the form $x_i \otimes x_j$ with $i \ne j$, and $\tilde{M}(\emptyset, (1))$ is spanned by the vectors $x_i \otimes x_i$.

\medskip

For the rest of the section we will recall some facts about permutation modules for $S_n$ and describe appropriate analogs for weighted permutation modules.

\subsubsection{Decomposition into irreducibles}

In the unweighted case, the decomposition of a permutation module into irreducibles given by the Kostka numbers which combinatorially count semistandard Young tableaux.  Explicitly we have $$M(\lambda) \cong \bigoplus_\mu K_{\lambda,\mu} S^\mu$$
we'll note however that this is just a special case of the Littlewood-Richardson rule (or even the Pieri rule) applied to the induction of the trivial representation of the Young subgroup corresponding to $\lambda$.

So unsurprisingly, the decomposition of weighted permutation modules will also be governed by Kostka numbers.  Our notation for polynomial representations will be more convenient for stating this result (and ultimately will be the case we care about), but we'll note that it holds for algebraic representations as well.

\begin{lemma}\label{multikostka}
$$\tilde{M}(\lambda^1, \lambda^2, \dots, \lambda^m) = \bigoplus_{\mu^1, \mu^2, \dots, \mu^m} K_{\mu^1,\lambda^1}K_{\mu^2,\lambda^2}\dots K_{\mu^m,\lambda^m} V^{\mu^1,\mu^2,\dots, \mu^m}$$
\end{lemma}

\noindent \textbf{Proof:} We've already seen that for algebraic $T \rtimes S_n$ modules induction is governed by the Littewood-Richardson rule on each $\mathbb{C}^\times$-weight separately.  As such, the decomposition of a weighted permutation modules is given by a product of Kostka numbers for each $\mathbb{C}^\times$ weight. \hfill $\square$

\subsubsection{$S_n$-invariants in weighted permutation modules}

The usual permutation module $M(\lambda)$ contains a one dimensional space of $S_n$-invariants. Indeed this is true anytime a representation of a finite group is constructed as the linearization of a transitive group action, and the space of invariants is spanned by the sum of the elements of the set being acted upon (or alternatively one can see this using Frobenius reciprocity as we are inducing up the trivial representation from a subgroup).

For weighted permutation modules, asking for a $T \rtimes S_n$ invariant vector is asking too much. Indeed this only occurs when $T$ acts trivially (which in the $\tilde{M}(\lambda, \mathbf{k})$ notation is when $\mathbf{k}$ is the zero vector), which just factors through the unweighted case. Nevertheless, we still have the following:

\begin{lemma}\label{perminv}
The permutation module $\tilde{M}(\lambda, \mathbf{k})$ has a one dimensional space of $S_n$-invariants.
\end{lemma}

\noindent \textbf{Proof:} By construction $M(\lambda, \mathbf{k})$ restricted to $S_n$ is isomorphic to $M(\lambda)$, which as we just said has a one dimensional space of invariants. \hfill $\square$

\subsubsection{Tensor products of weighted permutation modules}

Another key fact about permutation modules for $S_n$ is that the tensor product of two permutation modules is isomorphic to a direct sum of permutation modules.  More precisely, if $\lambda$ and $\mu$ are partitions of $n$ then $M(\lambda) \otimes M(\mu)$ decomposes into permutation modules as follows:

Define a \emph{tabloid} of type $(\lambda, \mu)$ to be a matrix of non-negative integers such that the $i$th row sums to $\lambda_i$ and the $j$th column sums to $\mu_j$ for all $i$ and $j$, and let $T(\lambda, \mu)$ be the set of all tabloids of type $(\lambda, \mu)$. Note that tabloids are just a convenient way of indexing the double cosets of the corresponding Young subgroups. For every tabloid $T \in T(\lambda, \mu)$, let $M(T)$ denote the corresponding permutation module where we think of the nonzero entries of $T$ as a composition of $n$. Then we have that
$$M(\lambda) \otimes M(\mu) = \bigoplus_{T(\lambda,\mu)} M(T)$$
and in fact this holds over arbitrary rings, although we'll stick to complex representations here (see \cite{Stanley} chapter 7). In particular we'll note that Lemma \ref{perminv} then says the space of $S_n$-invariants in $M(\lambda) \otimes M(\mu)$ is therefore equal to $|T(\lambda,\mu)|$, which is equal to the number of double cosets of the corresponding Young subgroups.

For weighted permutation modules essentially nothing changes, we just need to keep track of the weights.  If we are decomposing the tensor product $M(\lambda, \mathbf{k}) \otimes M(\lambda, \mathbf{k'})$, then define $T(\lambda, \mu)$ as before but this time define $M(T, \mathbf{k + k'})$ to be the weighted permutation module where if there is a nonzero entry in position $(i,j)$ we weight it by $k_i + k_j'$. We have the following lemma:

\begin{lemma}\label{permtensor}
$$M(\lambda, \mathbf{k}) \otimes M(\lambda, \mathbf{k'}) =  \bigoplus_{T(\lambda,\mu)} M(T, \mathbf{k + k'})$$
\end{lemma}

The proof is essentially identical to the unweighted case which is well known and a simple application of Mackey's theorem, so we'll omit it.

\medskip

\noindent \textbf{Example}: Suppose $n = 5$ and we want to decompose $M(\lambda,\mathbf{k}) \otimes M(\mu, \mathbf{k'})$ where $\lambda =  \mu = (3,2)$, $\mathbf{k} = (0,2)$ and $\mathbf{k'} = (1,4)$.  The relevant tabloids in $T(\lambda, \mu)$ are:

$$\begin{pmatrix}
3 & 0 \\
0 & 2
\end{pmatrix}
\hspace{1cm}
\begin{pmatrix}
2 & 1 \\
1 & 1
\end{pmatrix} 
\hspace{1cm}
\begin{pmatrix}
1 & 2 \\
2 & 0
\end{pmatrix}
$$
which give us the decomposition $$M(\lambda,\mathbf{k}) \otimes M(\mu, \mathbf{k'}) = M((3,2), (1,6)) \oplus M((2,1,1,1), (1,3,4,6)) \oplus M((1,2,2), (1,3,4))$$

\medskip

Next, suppose we want to decompose a tensor product of three or more permutation modules 
$$M(\lambda^1) \otimes M(\lambda^2) \otimes \dots \otimes M(\lambda^m)$$
as a direct sum of permutation modules.  Of course one could just repeatedly use the rule from above for a product of two permutation modules in terms of tabloids, however it turns out that it's often simpler to just do it in one step.

Define a \emph{multitabloid} of type $\lambda^1, \lambda^2, \dots, \lambda^m$ to be a $m$-dimensional array of nonnegative integers $b_{i_1,i_2,\dots, i_m}$  such that the $\ell$th generalized row sum in the $j$th direction is $\lambda^j_\ell$. That is,
$$\sum_{\substack{(i_1,i_2,\dots, i_m) \\ i_j = \ell}} b_{i_1,i_2,\dots, i_m} = \lambda^j_\ell$$
if we then let $T(\lambda^1, \lambda^2, \dots, \lambda^m)$ denote the set of all such multitabloids. Then by an easy induction on the number of terms $m$ we get the following corollary:

\begin{corollary} \ \\
\begin{enumerate}
\item For unweighted permutation modules $M(\lambda^1), M(\lambda^2), \dots , M(\lambda^m)$ we have: $$M(\lambda^1) \otimes M(\lambda^2) \otimes \dots \otimes M(\lambda^m) = \bigoplus_{B \in T(\lambda^1, \lambda^2, \dots, \lambda^m)} M(B)$$

\item For weighted permutation modules $M(\lambda^1, \mathbf{k}^1), M(\lambda^2, \mathbf{k}^2), \dots , M(\lambda^m,\mathbf{k}^m)$ we have $$M(\lambda^1, \mathbf{k}^1) \otimes M(\lambda^2, \mathbf{k}^2) \otimes \dots \otimes M(\lambda^m, \mathbf{k}^m) = \bigoplus_{B \in T(\lambda^1, \lambda^2, \dots, \lambda^m)} M(B, \mathbf{k}^1+ \mathbf{k}^2 + \dots + \mathbf{k}^m)$$

Where we view $T \in T(\lambda^1, \lambda^2, \dots, \lambda^m)$ as a composition by taking the nonzero entries, and $\mathbf{k}^1+ \mathbf{k}^2 + \dots + \mathbf{k}^m$ assigns the entry $b_{i_1,i_2,\dots, i_m}$ of $B$ the weight $k^1_{i_1} + k^2_{i_2}+ \dots +k^m_{i_m}$.

\end{enumerate}
\end{corollary}

\subsubsection{Stable tensor products for polynomial permutation modules}

We'll note our description of the decomposition of a tensor product of weighted permutation modules immediately gives us a compatible description for the polynomial permutation representations.  Indeed one just needs to recall that $\tilde{M}(\lambda^1, \lambda^2, \dots, \lambda^k)$ has an extra factor of size $(n-m)$ in degree zero when we write it as $M(\lambda, \mathbf{k})$. 

 The benefit of this notation for polynomial permutation representations is that it is independent of $n$ provided $n$ is large enough, and it is often easier to work with. For example, when we look at the matrices involved in decomposing a tensor product of permutation modules only thing that changes when we vary $n$ is the upper left entry. 
 
 For example recall from a previous example the decomposition
 
 $$V^{(1)} \otimes V^{(1)} = \tilde{M}(1) \otimes \tilde{M}(1)= \tilde{M}((1,1)) \oplus \tilde{M}(\emptyset, (1))$$
 in our other notation this is $$M((n-1,1), (0,1)) \otimes M((n-1,1), (0,1))$$
 and if $n$ is at least $2$ the corresponding tabloids are
 
 $$\begin{pmatrix}
n-2 & 1 \\
1 & 0
\end{pmatrix}
\hspace{1cm}
\begin{pmatrix}
n-1 & 0 \\
0 & 1
\end{pmatrix} $$
and we see that as $n$ changes the only thing that changes is the entry in the upper left corner.  Moreover we'll note that since the sum of the entries of these matrices sum to $n$, one can always recover that entry from knowing the rest of the matrix and $n$.

This motivates the notion of \emph{stable tabloids}, introduced in \cite{Harman} to study periodicity phenomena in the modular representation theory of symmetric groups, and independently in \cite{OZ1} to study certain symmetric functions related to $S_n$-representations.  If $\lambda = (\lambda_1, \lambda_2,\dots \lambda_m)$ and $\mu = (\mu_1, \mu_2, \dots ,\mu_\ell)$ are compositions of arbitrary sizes a stable tabloid of type $(\lambda,\mu)$ is a $(m+1)\times(\ell+1)$ matrix (with rows and columns indexed from zero to $m$ and $\ell$ respectively) such that:

\begin{enumerate}
\item The $(0,0)$ entry is empty and all other entries are non-negative integers.

\item For $i \ge 1$ the $i$th row sums to $\lambda_i$ and the $i$th column sums to $\mu_i$.

\end{enumerate}

The point being that for $n$ sufficiently large these are just obtained from deleting the upper left entries from the matrices in $T(\lambda[n],\mu[n])$ where $\lambda[n] = (n-|\lambda|, \lambda_1, \lambda_2,\dots \lambda_m)$ and $\mu[n] = (n-|\mu|, \mu_1, \mu_2, \dots ,\mu_\ell)$ are the padded partitions.

If we let $\tilde{T}(\lambda,\mu)$ denote the set of stable tabloids of type $(\lambda, \mu)$ then it follows (See \cite{Harman} or \cite{OZ1}) that the decomposition of the tensor product of (unweighted) permutation modules $M(\lambda[n]) \otimes M(\mu[n])$ for $n$ sufficiently large is just:

$$M(\lambda[n]) \otimes M(\mu[n]) = \bigoplus_{T \in \tilde{T}(\lambda,\mu)} M(T[n])$$
Where $T[n]$ is the the padded composition of $n$ obtained by taking the nonzero entries of $T$ along with an extra part of size $(n-|T|)$.

In the weighted case, Lemma \ref{permtensor} says that if we decompose a tensor product of two polynomial weighted permutation modules of the form $\tilde{M}(\lambda^1, \lambda^2, \dots, \lambda^m)$, then the decomposition into weighted permutation modules will again be indexed by weighted tabloids. To get the weights for a stable tabloid one just adds the corresponding weights on the rows and columns, with the caveat that the zeroth row and column have weight zero, in particular the ``extra" part of size $n-|T|$ always gets weight zero.

\medskip

\noindent \textbf{Remark:} We'll note that if $n$ is below the stable range we can still recover $T(\lambda[n],\mu[n])$ from $\tilde{T}(\lambda,\mu)$ just by throwing out those stable tabloids where the sum of the entries is larger than $n$ and then filling in the upper left corner appropriately.  Hence in this setting it will often be easy to deduce non-stable decompositions from the calculations in the stable range.

\medskip

As in the non-stable case, if we want to take a tensor product of three or more permutation modules
$$M(\lambda^1[n]) \otimes M(\lambda^2[n]) \otimes \dots \otimes M(\lambda^m[n])$$
 we can combine things into a single combinatorial object, a \emph{stable multitabloid}. If $\lambda_i$ has length $\ell_i$ these will be $m$-dimensional arrays of dimensions $(\ell_1 +1) \times (\ell_2 +1) \times \dots \times (\ell_m +1)$ such that:
 
 \begin{enumerate}
 \item The $(0,0,\dots, 0)$ position is left empty, and all other entries $b_{i_1,i_2, \dots i_m}$ are non-negative integers.
 
 \item The $r$th generalized row sum in the $j$th direction is $\lambda^j_r$. That is, $$\sum_ {\substack{(i_1,i_2, \dots i_m) \\ i_j = r}} b_{i_1,i_2, \dots i_m} = \lambda^j_r$$
 
 \end{enumerate}
 As before this can be seen by a simple inductive argument on the number of terms $m$, iterating the $m=2$ case at every step. Similarly to the case with two factors we can extend this to tensor products of weighted permutation modules by keeping track of the  weights accordingly.

\end{subsection}

\begin{subsection}{Symmetric powers}
We are now ready to start talking about symmetric powers of the defining representation $V$ and tensor products thereof.  The key observation is that as a representation of $T \rtimes S_n$ the symmetric power $Sym^k(V)$ decomposes as a direct sum of weighted permutation modules.

More precisely, let $A_k$ be the set of all $k$-tuples of non-negative integers $(a_1,a_2, \dots, a_k)$ such that 
$$a_1 + 2a_2 + 3a_3+\dots+ka_k = k$$
and let $A_k^n$ denote the subset of $A_k$ such that
$$a_1 + a_2 + a_3+\dots+a_k \le n$$
in particular note that if $n \ge k$ then $A_k^n = A_k$. The following proposition describes $Sym^k(V)$ as a sum of weighted permutation modules.

\begin{proposition}\label{symdec}
 $$Sym^k(V) \cong \bigoplus_{A_k^n} \tilde{M}((a_1),(a_2),(a_3), \dots, (a_k))$$

\end{proposition}

\noindent \textbf{Proof:} If $x_1,x_2,\dots x_n$ is the standard basis of $V$ then $Sym^k(V)$ is the space of homogeneous degree $k$ polynomials in the $x_i$'s.  Then for $(a_1,a_2, \dots, a_k) \in A_k^n$ consider the space spanned by all monomials $x_1^{b_1}x_2^{b_2} \dots x_n^{b_n}$ where $b_i = 1$ for $a_1$ values of $i$, $b_i = 2$ for $a_2$ values of $i$, and so on.  This space is isomorphic to the weighted permutation module $\tilde{M}((a_1),(a_2),(a_3), \dots, (a_k))$, with the monomials forming the permutation basis.  $Sym^k(V)$ is spanned by monomials and each monomial lies in some such space for a unique $(a_1,a_2, \dots, a_k) \in A_k^n$. \hfill $\square$

\medskip

We'll note that in this case the weighted permutation modules $\tilde{M}((a_1),(a_2),(a_3), \dots, (a_k))$ involved are all irreducible $T \rtimes S_n$ representations so we could also have written:

$$Sym^k(V) \cong \bigoplus_{A_k^n} V^{(a_1),(a_2),(a_3), \dots, (a_k)}$$

\subsubsection{Combinatorics of $A_k$}

In preparation for the next section we'd like to give two interpretations for $A_k$, which will generalize nicely to tensor products of symmetric powers.

\medskip

\noindent \textbf{Polynomial interpretation:} Given a sequence $(a_1,a_2, \dots, a_k)$ consider the polynomial $$a(x) := a_1x+a_2x^2+a_3x^3+\dots+a_kx^k$$
the condition that $(a_1,a_2, \dots, a_k) \in A_k$ just corresponds to the condition that $a'(1) = k$, and the condition that $(a_1,a_2, \dots, a_k) \in A_k^n$ adds the additional constraint that $a(1) \le n$. Hence $A_k^n$ can be naturally identified with the set of polynomials $a(x) \in \mathbb{Z}_{\ge 0}[x]$ such that $a(0)= 0$, $a(1) \le n$, and $a'(1) = k$.

\medskip

\noindent \textbf{Multiset partition interpretation:} Suppose $M$ is a multiset.  A \emph{multiset partition} of $M$ is a collection (i.e a multiset)  $\{ M_1, M_2, \dots, M_m\}$ of non-empty sub-multisets of $M$ such that the multiplicity of an element of $M$ is equal to the sum of the multiplicities of it in the $M_i$'s. 

$A_k$ can naturally be thought of as an indexing set for multiset partitions of the multiset $M = \{1,1,1, \dots, 1\} = \{1^k\}$ consisting of a single element with multiplicity $k$. Indeed an element $(a_1,a_2, \dots, a_k) \in A_k$ just corresponds to the unique multiset partition with $a_1$ parts of size $1$, $a_2$ parts of size 2, and so on. Under this identification $A_k^n$ is just indexing those multiset partitions of $M$ into at most $n$ parts.

\vspace{.5cm}

As mentioned before, in this case the decomposition into weighted permutation modules also gives the decomposition into irreducible $T \rtimes S_n$ modules. If we forget about the action of $T$ this gives us a combinatorial interpretation for the decomposition of $Sym^k(V)$ as a representation of $S_n$.

  For $(a_1, a_2, \dots, a_k) \in A_k^n$ define the \emph{associated composition of $n$} as $(n-(a_1 + a_2 + a_3+\dots+a_k), a_1, a_2, \dots, a_k)$, and define the associated composition of $n$ to a multiset partition of $M = \{1,1,1, \dots, 1\}$  via the identification above. 

\begin{corollary} If $\lambda$ is a partition of $n$ then the multiplicity of $S^\lambda$ in $Sym^k(V)$ is equal to the number of pairs $(P,T)$ where $P$ is a multiset partition of $\{1^k\}$ with at most $n$ parts and $T$ is a semistandard Young tableau of shape $\lambda$ and content equal to the associated composition of $n$ to $P$.

\end{corollary}

\noindent \textbf{Proof:} Proposition \ref{symdec} and the above interpretation tells us that $Sym^k(V)$ decomposes into weighted permutation modules indexed by these multiset partitions.  Restricting a weighted permutation modules to $S_n$ just gives an ordinary permutation module corresponding to the associated composition of $n$, and permutation modules decompose into irreducibles with multiplicities given by Kostka numbers which count semistandard Young tableaux.  \hfill $\square$

\medskip

In particular, the case where $\lambda = (n)$ and we are just looking at the space of symmetric group invariants is of particular interest so we'll state it as a separate corollary.

\begin{corollary}
The space of $S_n$-invariants in $Sym^k(V)$  is equal to the number of multiset partitions of $\{1^k\}$ with at most $n$ parts.

\end{corollary}

\end{subsection}

\begin{subsection}{Tensor products of symmetric powers}

Now let's extend this analysis from a single symmetric power of $V$ to a product of symmetric powers $$Sym^{k_1}(V) \otimes Sym^{k_2}(V) \otimes \dots \otimes Sym^{k_m}(V).$$
We already saw how to decompose a single symmetric power into weighted permutation modules, and we also saw how to decompose a tensor product of two weighted permutation modules, so all that remains is to put it all together and keep track of the terms to get a concise combinatorial description.

\subsubsection{The two factor case}

First let's go through how one would do this in the case where there are just two factors $Sym^{k_1}(V) \otimes Sym^{k_2}(V)$ just to demonstrate the idea.

\begin{enumerate}
\item First one would use Proposition \ref{symdec} to write  $$Sym^{k_i}(V) \cong \bigoplus_{A_{k_i}^n} \tilde{M}((a_1),(a_2),(a_3), \dots, (a_{k_i}))$$ for $i = 1,2$.

\item Then for each pair in $(a,a') \in A_{k_1}^n \times A_{k_2}^n$ we would decompose $$\tilde{M}((a_1),(a_2),(a_3), \dots, (a_{k_1})) \otimes \tilde{M}((a'_1),(a'_2),(a'_3), \dots, (a'_{k_2}))$$
according to Lemma \ref{permtensor} into weighted permutation modules indexed by stable tabloids of type $(a,a')$ (where we think of $a = (a_1,\dots, a_{k_1})$ and $a' = (a'_1,\dots, a'_{k_2})$ as compositions) with appropriate weights.

\end{enumerate}

Combining this into a single step by summing over all stable tabloids as $(a,a')$ varies over $A_{k_1}^n \times A_{k_2}^n$ we are summing over all $(k_1+1) \times (k_2+1)$ matrices with rows and columns indexed from $0$ to $k_1$ and $k_2$ respectively such that:

\begin{itemize}

\item The $(0,0)$ entry is blank, and all other entries are non-negative integers $b_{ij}$.

\item  The sum of the entries weighted by their row number is equal to $k_1$ and the sum of the entries weighted by their column number is equal to $k_2$. That is,

$$\sum_{i,j} i b_{ij} = k_1 \hspace{1cm} \sum_{i,j} j b_{ij} = k_2.$$

\item The sum of the entries of the matrix is at most $n$ (if we are in the stable range where $n > k_1 +k_2$ this condition is redundant).

\end{itemize}

For each such matrix we get a weighted permutation module $\tilde{M}(\lambda_1,\lambda_2, \dots, \lambda_{k_1 + k_2})$ where $\lambda_\ell$ is the partition obtained by taking the entries $b_{ij}$ with $i+j = \ell$.

\medskip

  For example if we want to decompose $Sym^2(V) \otimes Sym^2(V)$ there are $8$ such stable tabloids that appear:

   $$\begin{pmatrix}
\ & 2 & 0 \\
2 & 0 & 0 \\
0 & 0 & 0
\end{pmatrix}
\hspace{1cm}
\begin{pmatrix}
\ & 1 & 0 \\
1 & 1 & 0 \\
0 & 0 & 0
\end{pmatrix}
\hspace{1cm}
\begin{pmatrix}
\ & 0 & 0 \\
0 & 2 & 0 \\
0 & 0 & 0
\end{pmatrix}  
\hspace{1cm}
\begin{pmatrix}
\ & 1 & 0 \\
0 & 0 & 0 \\
0 & 1 & 0
\end{pmatrix}$$

   $$
\begin{pmatrix}
\ & 0 & 0 \\
1 & 0 & 1 \\
0 & 0 & 0
\end{pmatrix}
\hspace{1cm}
\begin{pmatrix}
\ & 2 & 0 \\
0 & 0 & 0 \\
1 & 0 & 0
\end{pmatrix} 
\hspace{1cm}
\begin{pmatrix}
\ & 0 & 1 \\
2 & 0 & 0 \\
0 & 0 & 0
\end{pmatrix}
\hspace{1cm}
\begin{pmatrix}
\ & 0 & 0 \\
0 & 0 & 0 \\
0 & 0 & 1
\end{pmatrix}  $$

We then read off the corresponding weighted permutation modules by looking at the diagonals that go from the lower left to the upper right to obtain the decomposition

\begin{equation}
\begin{split}
Sym^2(V) \otimes Sym^2(V) = \tilde{M}((2,2)) \oplus \tilde{M}((1,1),(1)) \oplus \tilde{M}(\emptyset, (2)) \\ \oplus \ 2\tilde{M}((1),\emptyset, (1)) \oplus 2\tilde{M}((2),(1)) \oplus \tilde{M}(\emptyset, \emptyset, \emptyset, (1))
\end{split}
\end{equation}
which holds for all $n \ge 4$. If $n$ is less than $4$ we obtain the decomposition by keeping just those terms where the sum of the sizes of the partitions involved is at most $n$.

\subsubsection{The general case}

Now let's extend this to the general case of $$Sym^{k_1}(V) \otimes Sym^{k_2}(V) \otimes \dots \otimes Sym^{k_m}(V).$$
Here the program will be basically identical to the two factor case

\begin{enumerate}
\item First one would use Proposition \ref{symdec} to write  $$Sym^{k_i}(V) \cong \bigoplus_{A_{k_i}^n} \tilde{M}((a_1),(a_2),(a_3), \dots, (a_{k_i}))$$ for $i = 1,2,\dots, m$.

\item Then for each $m$-tuple in $(a^1,a^2,\dots, a^m) \in A_{k_1}^n \times A_{k_2}^n\times \dots A_{k_m}^n$ we would decompose $$ \bigotimes_i \tilde{M}((a^i_1),(a^i_2),(a^i_3), \dots, (a^i_{k_i}))$$
according to Lemma \ref{permtensor} into weighted permutation modules indexed by stable multitabloids of type $(a^1,a^2, \dots, a^m)$ (where we think of $a^i = (a^i_1,\dots, a^i_{k_i})$ as compositions) with appropriate weights.

\end{enumerate}

As before we can combine this into a single step by counting up the multitabloids  as we vary our choice of $(a^1,a^2,\dots, a^m) \in A_{k_1}^n \times A_{k_2}^n\times \dots A_{k_m}^n$. Now the objects we care about are $(k_1+1)\times (k_2 +1) \times \dots \times (k_m+1)$ arrays with rows, columns, etc. indexed from $0$ to $k_1$, $k_2$, and so on respectively such that:

\begin{itemize}

\item The $(0,0, \dots , 0)$ entry is blank, and all other entries are non-negative integers $b_{i_1,i_2,\dots, i_m}$.

\item  The sum of the entries weighted by their $jth$ index value is equal to $k_j$  for all $j$. That is:

$$\sum_{i_1,i_2,\dots, i_m} i_j b_{i_1,i_2,\dots, i_m} = k_j$$

\item The sum of the entries of the array is at most $n$ (if we are in the stable range where $n > k_1 +k_2 + \dots +k_m$ this condition is redundant).

\end{itemize}

\medskip

If we let $A^n_{k_1,k_2,\dots,k_m}$ denote the set of such arrays then then putting everything together gives the following proposition describing the decomposition of a product of symmetric powers of $V$ into permutation modules.

\begin{proposition}\label{symproddec}
$$Sym^{k_1}(V) \otimes Sym^{k_2}(V) \otimes \dots \otimes Sym^{k_m}(V) = \bigoplus_{B \in A^n_{k_1,k_2,\dots,k_m}} \tilde{M}(\lambda_1,\lambda_2, \dots, \lambda_{k_1 + k_2+\dots +k_m})$$
where $\lambda_\ell = \lambda_\ell (B)$ is the partition obtained from the $B$ by taking the nonzero entries $b_{i_1,i_2,\dots,i_m}$ with $i_1+i_2+\dots+i_m = \ell$.

\end{proposition}

\subsubsection{Combinatorics of $A^n_{k_1,k_2,\dots,k_m}$}

As in the case of a single symmetric power, we will now give two interpretations for the set $A^n_{k_1,k_2,\dots,k_m}$, one in terms of polynomials and one in terms of multiset partitions.

\medskip

\noindent \textbf{Polynomial interpretation}:  Given an array $B$ of numbers $b_{i_1, i_2, \dots, i_m}$ (for $i_j \ge 0$) with finitely many nonzero entries we can construct the polynomial

$$P_B(x_1,x_2, \dots x_m) := \sum_{i_1, i_2, \dots, i_m} b_{i_1, i_2, \dots, i_m}x_1^{i_1}x_2^{i_2}\dots x_m^{i_m}$$
then the condition that $B \in A^n_{k_1,k_2,\dots,k_m}$ translates into

\begin{enumerate}
\item $P_B(x_1,x_2, \dots x_m)$ has constant term zero and all other coefficients are non-negative integers.

\item For $i =1, 2, \dots m$ $$ \frac{d}{dx_i}P_B(x_1,x_2, \dots x_m) |_{(x_1,x_2, \dots x_m) = (1,1,\dots 1)} = k_i$$

\item$ P_B(1,1,\dots 1) \le n$

\end{enumerate}

Hence we may identify $A^n_{k_1,k_2,\dots,k_m}$ with the set of such polynomials satisfying these three conditions.

\medskip

\noindent \textbf{Multiset partition interpretation:} $A^n_{k_1,k_2,\dots,k_m}$ has a natural bijection with the set of multiset partitions of $\{1^{k_1}, 2^{k_2}, \dots m^{k_m} \}$.  Explicitly, an array $B \in A^n_{k_1,k_2,\dots,k_m}$ with entries $b_{i_1, i_2, \dots, i_m}$ corresponds to the multiset partition where $\{1^{i_1}, 2^{i_2}, \dots, m^{i_m}\}$ appears $b_{i_1, i_2, \dots, i_m}$ times.

\vspace{.5cm}

We find the multiset partition interpretation to be more conceptually satisfying and will mostly use it to state results, however we've included the polynomial interpretation as we feel it may be more amenable to computations.

\medskip

 The \emph{type}, $\text{Type}(P)$, of a multiset partition $P \vDash M$ is the sequence $\lambda^1, \lambda^2, \lambda^3, \dots$ of partitions where $\text{Type}(P)^i = \lambda^i$ records the multiplicities of the parts of $P$ of size $i$. 
 
  For example, if $M = \{1,1,2\}$ and $P = \{\{1\},\{1\},\{2\}\}$ then $\text{Type}(P)$ is the sequence $\lambda^1, \lambda^2, \lambda^3, \dots$ where $\lambda^1 = (2,1)$ and all other partitions are the empty set. This is since $P$ has three parts of size $1$, two of which are the same as one another (i.e. $\{1\}$) and one that appears with multiplicity $1$ (i.e $\{2\}$). By convention we'll drop off the trailing empty sets and just write $\text{Type}(P) = ((2,1))$.  If instead we took $P' = \{1,1,2\}$ then $\text{Type}(P') = (\emptyset, \emptyset, (1))$.
 
Translating Proposition \ref{symproddec} into this language we get the following corollary:

\begin{corollary}\label{symtoperm}

$$Sym^{k_1}(V) \otimes Sym^{k_2}(V) \otimes \dots \otimes Sym^{k_m}(V) = \bigoplus_{\substack{P \vDash \{1^{k_1},2^{k_2}, \dots m^{k_m} \} \\ |P| \le n}}\tilde{M}(\text{Type}(P))$$

\end{corollary}

\medskip

We are now ready to give our main result, a combinatorial interpretation for the decomposition of a tensor product of symmetric powers into irreducible $T \rtimes S_n$ modules.

\begin{proposition}\label{maindecomp}
The multiplicity of $V^{\lambda^1, \lambda^2, \dots, \lambda^j}$  inside $$Sym^{k_1}(V) \otimes Sym^{k_2}(V) \otimes \dots \otimes Sym^{k_m}(V)$$
is equal to the number of tuples $(P,T_1,T_2, \dots T_j)$ such that $P$ is a multiset partition of $\{1^{k_1},2^{k_2}, \dots m^{k_m} \}$ with at most $n$ parts and $T_i$ is a semistandard Young tableau of shape $\lambda^i$ and content $\text{Type}(P)^i$.

\end{proposition}

\noindent \textbf{Proof:} Corollary \ref{symtoperm} gives us a decomposition of this tensor product of symmetric powers into permutation modules indexed by multiset partitions. Lemma \ref{multikostka} then gives us the decomposition of a permutation module in terms of Kostka numbers.  Finally, the combinatorial interpretation of Kostka numbers as counting semistandard Young tableaux gives the result. \hfill $\square$

\medskip

\noindent \textbf{Remark:} In the case where each $k_i$ is $1$ this refines Proposition $5.10$ in \cite{BHH} to the case of weighted permutation modules, but the proof is morally very similar.

\medskip

If we forget about the action of $T$, Corollary \ref{symtoperm} recovers an interpretation of the decomposition of $$Sym^{k_1}(V) \otimes Sym^{k_2}(V) \otimes \dots \otimes Sym^{k_m}(V)$$ into irreducible representations of $S_n$ due to Orellana and Zabrocki (\cite{OZ1} Theorem 5). 

 Define the \emph{unweighted type} $\text{UType}(P)$ of a multiset partition $P$ to be the partition recording the multiplicities of the multisets in $P$ (in other words just combine the partitions making up $\text{Type}(P)$ into a single partition). And let $\text{UType}(P)[n]$ be the corresponding partition of $n$ obtained by adding a part of size $(n-|P|)$ to $\text{UType}(P)$.

\begin{proposition}
The multiplicity of the Specht module $S^\lambda$  inside $$Sym^{k_1}(V) \otimes Sym^{k_2}(V) \otimes \dots \otimes Sym^{k_m}(V)$$
is equal to the number of tuples $(P,T)$ such that $P$ is a multiset partition of $\{1^{k_1},2^{k_2}, \dots m^{k_m} \}$ with at most $n$ parts and $T$ is a semistandard Young tableau of shape $\lambda$ and content $\text{UType}(P)[n]$.
\end{proposition}

\noindent \textbf{Proof:} Upon restricting to $S_n$  Corollary \ref{symtoperm} tells us that $$Sym^{k_1}(V) \otimes Sym^{k_2}(V) \otimes \dots \otimes Sym^{k_m}(V) = \bigoplus_{\substack{P \vDash \{1^{k_1},2^{k_2}, \dots m^{k_m} \} \\ |P| \le n}}M(\text{UType}(P)[n])$$
as a representation of $S_n$.  Then the decomposition of $M(\text{UType}(P)[n])$ is just given by appropriate Kostka numbers, which again count semistandard Young tableaux of fixed shape and content. \hfill $\square$

\medskip

\noindent \textbf{Remark:} In \cite{OZ1} Orellana and Zabrocki combine a pair $(P,T)$ into a single combinatorial object: a multiset tableau.  However morally these are the same description, and we'll leave it as an exercise to any interested parties to make explicit the bijection between these pairs and appropriate multiset tableaux.  

We'll also note that they stated their results in terms of a new basis of the ring of symmetric functions they defined corresponding to irreducible symmetric group representations of large symmetric groups.  In this context we'll mention that one should not expect such a basis to exist for $T \rtimes S_n$, as (among other reasons) it is possible for two elements of $T \rtimes S_n$ to be $GL_n$-conjugate but not $T \rtimes S_n$-conjugate.

\medskip

Finally, we'll close out this section by stating as a corollary the important special cases of the above propositions where we just look at the space of $S_n$-invariants.

\begin{corollary}

If $\mu = (\mu_1, \mu_2, \dots, \mu_\ell)$ is a partition of $k_1+k_2+\dots+k_m$ with at most $n$ parts the dimension of $S_n$-invariants in the symmetrized weight space

$$(Sym^{k_1}(V) \otimes Sym^{k_2}(V) \otimes \dots \otimes Sym^{k_m}(V))_{\bar{\mu}}$$
is equal to the number of multiset partitions of $\{1^{k_1},2^{k_2}, \dots m^{k_m} \}$ with parts of sizes $\mu_1, \mu_2, \dots, \mu_\ell$. In particular, the dimension of the space of $S_n$-invariants inside full space $$Sym^{k_1}(V) \otimes Sym^{k_2}(V) \otimes \dots \otimes Sym^{k_m}(V)$$
is equal to the number multiset partitions of $\{1^{k_1},2^{k_2}, \dots m^{k_m} \}$ with at most $n$ parts.

\end{corollary}

In particular if all the $k_i$'s are equal to $1$ and $m > n$ this recovers the fact that the space of $S_n$-invariants in $V^{\otimes m}$ is given by the Bell number $B_m$.  If we let the $k_i$'s be arbitrary but assume $n$ is sufficiently large this gives a representation theoretic context for the generalized Bell numbers studied in \cite{Griffiths}.

\end{subsection}
\end{section}

\begin{section}{Future directions and computations}\label{future}

The main motivating problem of finding a combinatorial interpretation for the restriction of an irreducible polynomial representation of $GL_n$ to $T \rtimes S_n$ remains open. In light of Proposition \ref{maindecomp} and the fact that weighted permutations satisfy an upper triangularity property with respect to the irreducibles it is reasonable to expect that such an interpretation could be formulated in terms of multiset partitions and Young tableaux, however for the time being such an interpretation remains elusive.

For the remainder of the paper we will outline some future directions of study beyond this motivating question. In particular we will outline a version of Schur-Weyl duality for $T \rtimes S_n$, highlight a connection between the restriction problem and the Foulkes conjecture, and include some low-degree computations of restriction from $GL_n$ to $T \rtimes S_n$.

\begin{subsection}{Schur-Weyl duality for $T \rtimes S_n$}

Now we'll breifly describe a version of Schur-Weyl duality for $T \rtimes S_n$.  As before let $V$ denote the defining representation of $GL_n$.  We have a tower of subgroups 
$$S_n \subset T\rtimes S_n \subset GL_n$$
acting on $V$ and hence on $V^{\otimes k}$.  Therefore we get a reverse inclusion of endomorphism algebras 
$$\End_{S_n}(V^{\otimes k}) \supset \End_{T \rtimes S_n}(V^{\otimes k}) \supset \End_{GL_n}(V^{\otimes k})$$
where the outer two terms are familiar instances of Schur-Weyl duality.

On the right we have classical Schur-Weyl duality in which $S_k$ acts $GL_n$-equivariantly on $V^{\otimes k}$ by permuting the factors.  For any $n$ and $k$ these maps span $ \End_{GL_n}(V^{\otimes k})$, and if $n \ge k$ this provides an isomorphism $\mathbb{C}[S_k] \cong  \End_{GL_n}(V^{\otimes k})$.

On the left the situation is similar. The partition algebra $Par(n,k)$ acts $S_n$-equivariantly on $V^{\otimes k}$ for all $n$ and $k$. The map $Par(n,k) \to \End_{S_n}(V^{\otimes k})$ is always surjective, moreover the kernel has an explicit combinatorial description. In particular, if $n \ge 2k$ this map is an isomorphism (See \cite{HR}, \cite{BH1}, \cite{BH2} for more details).

As such one may expect that we can explicitly describe the endomorphism algebra $\End_{T \rtimes S_n}(V^{\otimes k})$ at least for $n$ sufficiently large compared to $k$ and indeed this is the case.  Recall that $Par(n,k)$ has a basis indexed by set-partitions of the set $\{1,2,\dots,k\} \cup \{1',2',\dots k'\}$, we say that such a set-partition is \emph{balanced} if each part $P$ of the partition satisfies $|P \cap \{1,2,\dots,k\}| = |P \cap \{1',2',\dots,k'\}|$. 

\medskip

The following proposition describes a version of Schur-Weyl duality for $T\rtimes S_n$:
 
\begin{proposition} Let $Par^{bal}(k)$ denote the subspace of $Par(n,k)$ spanned by the balanced set partitions.  The following holds:
\begin{enumerate}
\item $Par^{bal}(k)$ is a subalgebra of $Par(n,k)$.
\item The obvious identification of these subspaces for $Par(n,k)$ and $Par(m,k)$ is an isomorphism of algebras (motivating the notation $Par^{bal}(k)$ not involving $n$).
\item  $Par^{bal}(k)$ acts $T\rtimes S_n$-equivariantly on $V^{\otimes k}$ giving an algebra homomorphism $Par^{bal}(k) \to  \End_{T \rtimes S_n}(V^{\otimes k})$, this is surjective for all values of $n$ and $k$ and is an isomorphism whenever $n \ge k$.

\end{enumerate}
\end{proposition}

\noindent \textbf{Sketch of proof:} For parts $1$ and $3$ we use the description of Schur-Weyl duality for $S_n$ and check that $Par^{bal}(k)$ is exactly the subalgebra of $Par(n,k)$ that preserves $S_n$-orbits of $T$-weights. For part $2$ one just notes that in the partition algebra $Par(n,k)$, the dependence on $n$ arises when we stack two partition diagrams and there are isolated components in the middle, but for balanced partitions that can never happen. \hfill $\square$

\medskip

We didn't pursue this direction any further, but here are a few problems we think may be interesting:

\medskip

\noindent \textbf{Problem:} Describe $Par^{bal}(k)$ by generators and relations.

\medskip

\noindent \textbf{Problem:} Is there a double centralizer property? In other words, does the image of $T \rtimes S_n$ span the algebra of $Par^{bal}(k)$-endomorphisms of $V^{\otimes k}$? If so, use this to describe the representations of $Par^{bal}(k)$.

\medskip

\noindent \textbf{Problem:} Describe the kernel of the natural map $Par^{bal}(k) \to  \End_{T \rtimes S_n}(V^{\otimes k})$ for $n < k$.

\medskip

\noindent \textbf{Problem:} Describe the $T \rtimes S_n$ endomorphism algebras of $V^{\otimes r} \otimes V^{* \otimes s}$. We'll note that this should be something in between a partition algebra and a walled Brauer algebra.
\end{subsection}

\medskip

\noindent \textbf{Problem:} Describe diagrammatically the endomorphism algebras of $$Sym^{k_1}(V) \otimes Sym^{k_2}(V) \otimes \dots \otimes Sym^{k_m}(V)$$ as $S_n$ and $T \rtimes S_n$ representations.  We'll go ahead and coin the terms ``multiset partition algebras" and ``balanced multiset partition algebras".

\begin{subsection}{$S_n$-invariants, plethysm, and Foulkes conjecture}

A particularly important subcase of general problem of restricting an irreducible representation $W(\lambda)$ from $GL_n$ to $T \rtimes S_n$ is to understand the the space of $S_n$-invariants in each symmetrized weight space $W(\lambda)_{\bar{\mu}}$ as $\mu$ varies over $S_n$-orbits of weights.

If $\mu$ has $m_1$ parts of size $1$, $m_2$ parts of size $2$, and so on.  Let 

$$S_{\mu} := S_1^{m_1} \times S_2^{m_2} \times \dots \times S_k^{m_k}$$
denote the Young subgroup of $S_{|\mu|}$ corresponding to $\mu$, and let $N(\mu)$ denote its normalizer in $S_{|\mu|}$.  Note that $N(\mu)$ is just the product of the wreath products $S_i \wr S_{m_i}$. 

The following proposition, essentially due to Gay, relates the space of $S_n$-invariants in the symmetrized weight space $W(\lambda)_{\bar{\mu}}$to the representation theory of $S_{|\lambda|}$.

\begin{proposition} \textbf{(\cite{Gay} Theorem 2)} The dimension of the space of $S_n$-invariants in $W(\lambda)_{\bar{\mu}}$is equal to the multiplicity of the Specht module $S^\lambda$ in $\text{Ind}_{N(\mu)}^{S_{|\mu|}}(\mathbf{1})$.
\end{proposition}

\noindent \textbf{Remark:} We'll note that Gay only considered the ``zero" weight space case where $\mu = (a,a,a,\dots,a)$.  However his proof easily generalizes, and moreover one can reduce the general case to this by a straightforward application of the Littlewood-Richardson rule for restricting to the Levi subgroup $GL_{m_1} \times GL_{m_2} \times \dots \times GL_{m_k}$.

\medskip

Finding a combinatorial interpretation of such modules induced from normalizers of Young subgroups is an important open problem in combinatorial representation theory which in the language of symmetric functions is equivalent to decomposing the plethysms $h_a[h_b]$ into Schur functions. 

In particular one important conjecture in this area is Foulkes conjecture from 1950, about embedding one such induced module into another.  Translated into the language of this paper Foulkes conjecture is the following:

\begin{conjecture} \textbf{(Foulkes conjecture \cite{Foulkes})} If $a < b$ then the space of $S_n$-invariants in the symmetrized weight space of weight $(a,a,\dots, a, 0,0,\dots,0)$ is at least as large as the space of $S_n$-invariants in the symmetrized weight space of weight $(b,b,\dots, b, 0,0,\dots,0)$ for any polynomial representation of $GL_n$ of degree $ab$.

\end{conjecture}

We'll note that the usual representation theoretic formulations of this conjecture are either entirely about general linear group representations or entirely about symmetric group representations. As far as we can tell this ``mixed"  formulation of the conjecture seems to be missing from much of the literature on Foulkes conjecture.

  It suggests a possible approach to the conjecture by studying the representation theory of the spherical subalgebra $eAe$, where $$A = U(\mathfrak{gl}_n) \rtimes \mathbb{C}(S_n) \hspace{.5cm} \text{ and }  \hspace{.5cm} e = \frac{1}{n!}\sum_{\sigma \in S_n} \sigma$$
which naturally acts on the space of $S_n$-invariants inside a $GL_n$ representation. 

\medskip

Using Proposition \ref{maindecomp} we can give a purely combinatorial formulation of a weak version of Foulkes conjecture for tensor products of symmetric powers:

\begin{conjecture} \textbf{(Weak Foulkes conjecture)} If $a < b$ and $M$ is a multiset of size $ab$ then the number of multiset partitions of $M$ into $b$ parts of size $a$ at least as large as the number of multiset partitions of $M$ into $a$ parts of size $b$.
\end{conjecture}

Note that in terms of symmetric functions this is equivalent to the conjecture that $h_b[h_a]-h_b[h_a]$ has nonnegative coefficients when expressed in the basis of monomial symmetric functions (whereas the full Foulkes conjecture says it has non-negative coefficients in the basis of Schur functions).

\end{subsection}

\subsection{Low degree calculations}
Proposition \ref{maindecomp} combined with the Jacobi-Trudi identity gives us a combinatorial method for computing the decomposition of an irreducible polynomial representation of $GL_n$ to $T \rtimes S_n$ (although not a positive combinatorial formula). Here are some explicit calculations for polynomial representations of degree at most $4$, where we always assume we are in the stable range with $n$ larger than the degree.

$$ W( \emptyset) \longrightarrow V^\emptyset $$
$$W(1) \longrightarrow V^{(1)} $$
$$W(2) \longrightarrow V^{(2)} \oplus V^{\emptyset, (1)}$$
$$ W(1,1) \longrightarrow V^{(1,1)}$$
$$W(3) \longrightarrow V^{(3)} \oplus V^{(1),(1)} \oplus V^{\emptyset,\emptyset, (1)}$$
$$W(2,1) \longrightarrow V^{(2,1)} \oplus V^{(1),(1)}$$
$$W(1,1,1) \longrightarrow V^{(1,1,1)}$$
$$W(4) \longrightarrow V^{(4)} \oplus V^{(2),(1)} \oplus V^{\emptyset, (2)} \oplus V^{(1), \emptyset, (1)} \oplus V^{\emptyset,\emptyset, \emptyset, (1)}$$
$$W(3,1) \longrightarrow V^{(3,1)} \oplus V^{(1,1),(1)} \oplus V^{(2),(1)} \oplus V^{\emptyset, (1,1)} \oplus V^{(1),\emptyset,(1)} $$
$$ W(2,2) \longrightarrow V^{(2,2)} \oplus V^{(2),(1)} \oplus V^{\emptyset, (2)}$$
$$W(2,1,1) \longrightarrow V^{(2,1,1)} \oplus V^{(1,1),(1)}$$
$$W(1,1,1,1) \longrightarrow V^{(1,1,1,1)}$$

\end{section}

\end{document}